\magnification=\magstep1
\documentstyle{amsppt}

\topmatter
\title 
Deformations of canonical singularities
\endtitle
\author Yujiro Kawamata
\endauthor

\rightheadtext{Singularities}

\address Department of Mathematical Sciences, University of Tokyo, 
Komaba, Meguro, Tokyo, 153, Japan \endaddress
\email kawamata\@ms.u-tokyo.ac.jp\endemail


\abstract
We prove that small deformations of canonical singularities are canonical.
\endabstract

\endtopmatter

\document

\head 1. Introduction
\endhead

M. Reid defined the concept of canonical singularities 
as those which appear on canonical models
of algebraic varieties whose canonical rings are finitely generated.
From the view point of the classification theory of algebraic varieties,
it is expected that small deformations of a variety which has only canonical
singularities have only canonical singularities.
The main result of this paper states that this is the case:

\proclaim{Main Theorem}
Let $\pi: X \to S$ be a flat morphism from a germ of a variety 
$(X, x_0)$ to a germ of a 
smooth curve $(S, s_0)$ whose special fiber $X_0 = \pi^{-1}(s_0)$
has only canonical singularities.
Then $X$ has only canonical singularities.
In particular, the fibers $X_s = \pi^{-1}(s)$ 
have only canonical singularities.
\endproclaim

A normal variety $X$ is said to have only canonical singularities if
the canonical divisor $K_X$ is a $\Bbb Q$-Cartier divisor 
and if the discrepancy divisor
$K_Y - \mu^*K_X$ is effective for a resolution of singularities
$\mu: Y \to X$ (cf. [KMM]).
For example, 2-dimensional canonical singularities are nothing 
but Du Val (or rational double) singularities. 
It is well known that deformations of Du Val singularities are 
Du Val singularities.
Our theorem was not known even for the 3-dimensional case.
Similar problem for terminal singularities is still open.
A result is already obtained for log terminal singularities ([I]).

The proof of the main theorem follows closely the argument of Siu [Si] who 
proved the following: if $\pi: X \to S$ is a smooth projective family of
varieties over a base space such that all the fibers $X_s$ are of general type,
then the plurigenus $P_m(X_s)$ of the fibers is constant for each 
positive integer $m$.
The difference is that our method is algebraic instead of analytic.
We can also adapt the analytic proof of Siu's theorem to the algebraic 
one by using our method (Theorem 6).

The concept of multiplier ideal sheaf introduced by Nadel in the
analytic setting plays an important role in the proof.
The point is that the multiplier ideal sheaf has a canonical
embedding to the structure sheaf, and one can consider it
independently from the original linear system.
In order to deal with an infinite sequence of linear systems, 
we consider an increasing sequence of multiplier ideal sheaves
instead of considering the infinite sum of metrics as in [Si].

\head 2. Proof of the Main Theorem
\endhead

Let $\mu: Y \to X$ be a projective birational morphism from a smooth variety 
such that $Y_0 = \mu^*X_0 = X'_0 + E$ 
is a normal crossing divisor, where $X'_0$ is the strict transform of
$X_0$. 
We denote $\pi' = \pi \circ \mu$.
We set $\mu^*K_{X/S} = K_{Y/S} + \Gamma$.
We have to prove that $- \Gamma$ is effective.

It is easy to see that $X$ is normal.
The canonical divisor $K_X$ is known to be $\Bbb Q$-Cartier ([St, \S 6]).
Thus what we have to prove is that any section of $\Cal O_X(mK_{X/S})$ for
any positive integer $m$ can be 
lifted to a section of $\Cal O_Y(mK_{Y/S})$.
Since the sections of $\Cal O_{X_0}(mK_{X_0})$ can be lifted to that of 
$\Cal O_{X'_0}(mK_{X'_0})$ by the assumption,
the problem is reduced to extend sections of $\Cal O_{X'_0}(mK_{X'_0})$ to
those of $\Cal O_Y(mK_{Y/S})$.

For positive integers $m$, 
we construct a series of projective birational morphisms $f_m: Y_m \to Y$ from 
smooth varieties which satisfy the following conditions (1) through (11).
We denote $f'_m = \mu \circ f_m$ and $\pi_m = \pi \circ f'_m$. 
A divisor $\Delta_m$ on $Y_m$ is defined by
$K_{Y_m/S} + \Delta_m = f_m^*K_{Y/S}$.
We note that $- \Delta_m$ is effective.

(1) There is a decomposition 
$$
f_m^*(mK_{Y/S}) = P_m + M_m
$$ 
in the group of divisors $\text{Div}(Y_m)$.

(2) $P_m$ is $f'_m$-free, i.e., the natural homomorphism
$f_m^{\prime *}f'_{m*}\Cal O_{Y_m}(P_m) \to \Cal O_{Y_m}(P_m)$ 
is surjective.

(3) $M_m$ is an effective divisor whose support is contained in a 
normal crossing divisor $F_m$.

(4) $M_m$ is the $f_m'$-fixed part of $mK_{Y/S}$, that is, 
the natural homomorphism 
$f'_{m*}\Cal O_{Y_m}(P_m) \to \mu_*\Cal O_Y(mK_{Y/S})$ 
is bijective.

(5) The supports of the divisors $f_m^*Y_0$, $f_m^*\Gamma$ and 
$\Delta_m$ are contained in $F_m$.

For any very $\mu$-ample divisor $A$ on $Y$, 
there exist a positive integer $m_0$ and an effective divisor 
$B$ whose support does not contain $X'_0$ 
such that $m_0K_{Y/S} \sim A + B$, because $X$ is a germ.
By definition, we have $f_{m_0}^*B \ge M_{m_0}$.

(6) The support of the divisor $f_m^*B$ is contained in $F_m$.

We also assume that there is a projective birational morphism 
$g_{km}: Y_m \to Y_k$ such that 
$f_m = f_k \circ g_{km}$ for each pair of integers $k,m$ with $k < m$.

(7) The supports of $g_{km}^*M_k$ and $g_{km}^*F_k$ are contained in $F_m$.

We have $g_{km}^*M_k/k \ge M_m/m$ if $k \vert m$.
Let $Y_{m,0} = f_m^*Y_0 = X'_{m,0} + E_m$, where 
$X'_{m,0}$ is the strict transfrom of $X'_0$.
Let $Y'_{m,0} = f_m^*X'_0$.
Let $\mu_0: X'_0 \to X_0$,
$f_{m,0}: X'_{m,0} \to X'_0$, $f'_{m,0}: X'_{m,0} \to X_0$,
$\bar f_{m,0}: Y'_{m,0} \to X'_0$ 
$\bar f'_{m,0}: Y'_{m,0} \to X_0$, and
$g_{km,0}: X'_{m,0} \to X'_{k,0}$
be the induced morphisms:
$$
\align
&Y_m @>{g_{km}}>> Y_k @>{f_k}>> Y @>{\mu}>> X @>{\pi}>> S \\
&X'_{m,0} @>{g_{km,0}}>> X'_{k,0} @>{f_{k,0}}>> X'_0 @>{\mu_0}>> X_0
\endalign
$$

We define $\Gamma_0$ and $\Delta_{m,0}$ by
$K_{X'_0} + \Gamma_0 = \mu_0^*K_{X_0}$ and
$K_{X'_{m,0}} + \Delta_{m,0} = f_{m,0}^*K_{X'_0}$.
Since $X_0$ has only canonical singularities, 
we have $\Gamma_0 \le 0$. 
We also have $\Delta_{m,0} \le 0$.
By the adjunction formula, we have
$\Gamma_0 = (\Gamma + E) \vert_{X'_0}$ and
$\Delta_{m,0} = (\Delta_m + E_m - f_m^*E) \vert_{X'_{m,0}}$.
Let $B_0 = (B - m_0E)\vert_{X'_0}$ and
$B_{m,0} = f_m^*(B - m_0E)\vert_{X'_{m,0}}$.
Let $F_{m,0} = (F_m - X'_{m,0}) \cap X'_{m,0}$.

(8) There is another decomposition 
$$
f_{m,0}^*(mK_{X'_0}) = Q_m + N_m
$$ 
in $\text{Div}(X'_{m,0})$.

(9) $Q_m$ is $f'_{m,0}$-free, i.e., the natural homomorphism
$f_{m,0}^{\prime *}f'_{m,0*}\Cal O_{X'_{m,0}}(Q_m) \to \Cal O_{X'_{m,0}}(Q_m)$ 
is surjective.

(10) $N_m$ is an effective divisor whose support is contained in the
normal crossing divisor $F_{m,0}$.

(11) $N_m$ is the $f'_{m,0}$-fixed part of $mK_{X'_0}$, that is, 
the natural homomorphism 
$f'_{m,0*}\Cal O_{X'_{m,0}}(Q_m) \to \mu_{0*}\Cal O_{X'_0}(mK_{X'_0})$ 
is bijective.

There is an injective homomorphism 
$\mu_*\Cal O_Y(mK_{Y/S}) \to \Cal O_X(mK_{X/S})$,
because $\Cal O_X(mK_{X/S})$ is a reflexive sheaf.
Since $X_0$ has ony canonical singularities, 
we have a natural homomorphism 
$\Cal O_X(mK_{X/S}) \to \mu_{0*}\Cal O_{X'_0}(mK_{X'_0})$
for $m \ge 0$.
Therefore, there is a natural homomorphism
$$
\mu_*\Cal O_Y(mK_{Y/S}) \to \mu_{0*}\Cal O_{X'_0}(mK_{X'_0}).
$$
Thus we have $P_m \vert_{X'_{m,0}} \le Q_m$.
We set $M_{m,0} = (M_m - mf_m^*E) \vert_{X'_{m,0}}$.
Since $(K_{Y/S} - E)\vert_{X'_0} = K_{X'_0}$, 
we have $M_{m,0}/m \ge N_m/m$.
We have $B_{m,0}/m_0 \ge M_{m,0}/m$ if $m_0 \vert m$.

\definition{Definition 1}
Let $f: Y \to X$ be a projective birational morphism of smooth varieties,
and let $\pi: X \to S$ be another projective morphism
to a base space $S$, which is assumed to be a quasi-projective variety or 
its germ at a point.
Let $f' = \pi \circ f$.
Let $\Delta$ be a divisor on $Y$ defined by
$K_Y + \Delta = f^*K_X$.
Let $D$ be a divisor on $X$.
Assume that $f^*D = P + N$, where 
$P$ is a $f'$-nef and $f'$-big $\Bbb R$-divisor, 
and $N$ is an effective $\Bbb R$-divisor 
whose support is a normal crossing divisor.
The {\it multiplier ideal sheaf} $\Cal I_N$ on $X$ is defined by the following 
formula: $\Cal I_N = f_*\Cal O_Y(\ulcorner - N - \Delta \urcorner)$.
Since the effective part of the divisor $\ulcorner - N - \Delta \urcorner$
is exceptional for $f$, $\Cal I_N$ is an ideal sheaf of $\Cal O_X$ 
which is determined by $X$ and $N$.
Alternatively, we have
$$
f_*\Cal O_Y(\ulcorner P \urcorner + K_Y) = \Cal I_N(D + K_X).
$$

If $g: Y' \to Y$ is a projective birational morphism from another smooth 
variety such that the support of $N' = g^*N$ is a normal crossing divisor,
then we have $\Cal I_{N'} = \Cal I_N$,
because $\ulcorner g^*P \urcorner + K_{Y'} \ge 
g^*(\ulcorner P \urcorner + K_Y)$.

Since $R^pf'_*\Cal O_Y(\ulcorner P \urcorner + K_Y) = 0$ and
$R^pf_*\Cal O_Y(\ulcorner P \urcorner + K_Y) = 0$ for $p > 0$, we have
$R^p\pi_*(\Cal I_N(D + K_X)) = 0$ for $p > 0$.
\enddefinition

\proclaim{Theorem 2}
Let $\pi: X \to S$ be a projective morphism from a smooth variety to
a base space.
Then there exists a $\pi$-ample divisor $A$ which depends only on $\pi$ 
such that the natural homomorphism
$\pi^*\pi_*(\Cal I_N(A + D + K_X)) \to \Cal I_N(A + D + K_X)$ is
surjective for any $f: Y \to X$, $D$ and $N$ as in Definition 1.
\endproclaim

\demo{Proof}
Let $x$ be any point of $X$, and $\frak m_x$ the corresponding sheaf of ideals
of $\Cal O_X$.
Let $g: X' \to X$ be the blowing up at the ideal sheaf
$\frak m_x$.
Let $E$ be the coresponding Cartier divisor on $X'$; we have
$\frak m_x\Cal O_{X'} = \Cal O_{X'}(-E)$.
We note that $- E$ is $g$-ample.
Let $g' = \pi \circ g$.

We take a $\pi$-ample divisor $A$ on $X$ such that the natural 
homomorphism
$$
g^{\prime *}g'_*(\Cal O_{X'}(- E + g^*A)) \to \Cal O_{X'}(- E + g^*A)
$$ 
is surjective. 
Moreover, we assume that there is a surjective homomorphism
$\Cal O_S^{\oplus \ell} \to g'_*\Cal O_{X'}(- E + g^*A)$.
Let $\Cal F$ be the kernel of the induced homomorphism 
$\Cal O_{X'}^{\oplus \ell} \to \Cal O_{X'}(- E + g^*A)$,
which is a locally free sheaf on $X'$.
Let $V = \Bbb P(\Cal F)$ be the associated projective space
bundle with the projection $\sigma: V\to X'$ and
the tautological sheaf $\Cal O_V(H)$.

We may assume that there exists a projective birational morphism 
$h: Y \to X'$ such that $f = g \circ h$.
We set $W = \Bbb P(h^*\Cal F)$.
Let $\bar h: W \to V$ and $\tau: W \to Y$ be natural morphisms:
$$
\CD
W @>{\bar h}>> V \\
@V{\tau}VV       @V{\sigma}VV \\
Y @>{h}>>      X' @>g>> X @>{\pi}>> S.
\endCD
$$

Let $\Cal J$ be an ideal shaef of $\Cal O_{X'}$ defined by 
$\Cal J = h_*\Cal O_Y(\ulcorner - N - \Delta \urcorner)$.
Then we have $h_*\Cal O_Y(\ulcorner P \urcorner + K_Y) = \Cal J(g^*(D+K_X))$ 
and $g_*\Cal J = \Cal I_N$.
We take a positive integer $b$ such that
the divisor $H - b\sigma^*E - K_{V/X'}$ on $V$ is $(g \circ \sigma)$-ample.
We also assume now that the divisor $- (b+1)E + g^*A$ is $g'$-ample.
Then we have for $p > 0$
$$
\align
&R^pf_*(h^*\Cal F(-bh^*E + \ulcorner P \urcorner + K_Y)) \\
&= R^p(f \circ \tau)_*\Cal O_W(\bar h^*(H - b\sigma^*E - K_{V/X'})
+ \ulcorner \tau^*P \urcorner + K_W) = 0 \\
&R^ph_*(h^*\Cal F(-bh^*E + \ulcorner P \urcorner + K_Y)) = 0
\endalign
$$
hence
$$
R^1g_*(\Cal J \otimes \Cal F(-bE)) = 0.
$$
Similarly, we have 
$$
R^1g'_*(\Cal J(-(b+1)E + g^*(A+D+K_X))) 
= R^1g_*(\Cal J(-(b+1)E)) 
= 0.
$$ 
Therefore, from an exact sequence
$$
0 \to \Cal J \otimes \Cal F(-bE) \to \Cal J^{\oplus \ell}(-bE)
\to \Cal J(- (b+1)E + g^*A) \to 0,
$$
we obtain an exact sequence
$$
0 \to g_*(\Cal J \otimes \Cal F(-bE)) \to 
g_*(\Cal J^{\oplus \ell}(-bE))
\to g_*(\Cal J(- (b+1)E+g^*A)) \to 0.
$$
Thus 
$$
g_*(\Cal J(- (b+1)E)) = g_*(\Cal J(- bE)) \cdot \frak m_x \subset 
\Cal I_N \cdot \frak m_x.
$$

On the other hand, we have surjective homomorphisms
$$
\align
&g'_*(\Cal J(g^*(A+D+K_X))) 
\to g'_*(\Cal J(g^*(A+D+K_X)) \otimes \Cal O_{(b+1)E}) \\
&g_*(\Cal J(g^*(A+D+K_X))) 
\to g_*(\Cal J(g^*(A+D+K_X)) \otimes \Cal O_{(b+1)E})
\endalign
$$
whose images are naturally identified.
Hence we obtain natural surjective homomorphisms
$$
\align
&\pi_*(\Cal I_N(A+D+K_X)) \\
&\to g_*(\Cal J(g^*(A+D+K_X)))/
g_*(\Cal J(- (b+1)E + g^*(A+D+K_X))) \\
&\to I_N(A+D+K_X) \otimes \Cal O_X/\frak m_x.
\endalign
$$
\hfill $\square$
\enddemo

\definition{Remark 3}
The above proof can be easily extended to the case in which $X$ is singular.
As long as $X$ is smooth, 
one can use the Koszul resolution of $\Cal O_{X'}(-E)$
instead of the sheaf $\Cal F$ of the above proof 
in order to obtain a more precise condition on the divisor $A$.
\enddefinition

Since $K_{Y/S} \vert_{X'_0} = K_{X'_0} + E \vert_{X'_0}$, we have
$$
\align
f_{m*}\Cal O_{Y_m}(\ulcorner kP_m/m \urcorner + K_{Y_m/S})
&= \Cal I_{kM_m/m}((k+1)K_{Y/S}) \\
f_{m,0*}\Cal O_{X'_{m,0}}(\ulcorner kP_m/m \vert_{X'_{m,0}} \urcorner 
+ K_{X'_{m,0}})
&= \Cal I_{kM_{m,0}/m}((k+1)K_{X'_0}) \\
f_{m,0*}\Cal O_{X'_{m,0}}(\ulcorner kQ_m/m \urcorner 
+ K_{X'_{m,0}})
&= \Cal I_{kN_m/m}((k+1)K_{X'_0}).
\endalign
$$
We note that the support of $\Cal O_Y/\Cal I_{kM_m/m}$ does not 
contain $X'_0$.

Since $g_{mm',0}^*N_m/m \ge N_{m'}/m'$ if $m \vert m'$, 
we have $\Cal I_{kN_m/m} \subset \Cal I_{kN_{m'}/m'}$ 
for any positive integer $k$.
We define
$$
\Cal I_{kN} = \bigcup_{m=1}^{\infty} \Cal I_{kN_m/m}.
$$ 
Since $\Cal O_{X'_0}$ is noetherian, there exists a positive integer $m$,
which depends on $k$, such that $\Cal I_{kN} = \Cal I_{kN_m/m}$.
Similarly, we define 
$\Cal I_{kM} = \bigcup_{m=1}^{\infty} \Cal I_{kM_{m,0}/m}$ and
$\Cal I_{k\bar M} = \bigcup_{m=1}^{\infty} \Cal I_{kM_m/m}$.

\proclaim{Lemma 4}
$$
\align
&\text{Im}(\mu_*\Cal O_Y((k+1)K_{Y/S}) \to 
\mu_{0*}\Cal O_{X'_0}((k+1)K_{X'_0})) \\
&= \mu_{0*}(\Cal I_{kM}((k+1)K_{X'_0})) \\
&= \mu_{0*}(\Cal I_{(k+1)M}((k+1)K_{X'_0})).
\endalign
$$
\endproclaim

\demo{Proof}
Since $\Cal I_{(k+1)M} \subset \Cal I_{kM}$, we have 
$$
\mu_{0*}(\Cal I_{(k+1)M}((k+1)K_{X'_0})) \subset 
\mu_{0*}(\Cal I_{kM}((k+1)K_{X'_0})).
$$

If $s \in \mu_{0*}\Cal O_{X'_0}((k+1)K_{X'_0})$
is the image of $\bar s \in \mu_*\Cal O_Y((k+1)K_{Y/S})$, then
$\text{div}(f_{k+1}^*\bar s) \ge M_{k+1}$ on $Y_{k+1}$.
Since $K_{Y/S} \vert_{X'_0} = K_{X'_0} + E\vert_{X'_0}$
and $M_{k+1,0} = (M_{k+1} - (k+1)f_{k+1}^*E)\vert_{X'_{k+1,0}}$, 
we have $\text{div}(f_{k+1,0}^*s) \ge M_{k+1,0}$,
hence $s \in \mu_{0*}I_{M_{k+1,0}}((k+1)K_{X_0})$.
Since $I_{M_{k+1,0}} \subset I_{(k+1)M}$, we have 
$$
\text{Im}(\mu_*\Cal O_Y((k+1)K_{Y/S}) \to 
\mu_{0*}\Cal O_{X'_0}((k+1)K_{X'_0})) 
\subset \mu_{0*}(\Cal I_{(k+1)M}((k+1)K_{X'_0})).
$$

We shall prove that
$$
\mu_{0*}(\Cal I_{kM}((k+1)K_{X'_0})) \subset
\text{Im}(\mu_*\Cal O_Y((k+1)K_{Y/S}) \to 
\mu_{0*}\Cal O_{X'_0}((k+1)K_{X'_0})).
$$
We take a positive integer $m$ such that 
$\Cal I_{kM} = \Cal I_{kM_{m,0}/m}$ and
$\Cal I_{k\bar M} = \Cal I_{kM_m/m}$.
Since $R^1f_{m*}\Cal O_{Y_m}(\ulcorner kP_m/m \urcorner + K_{Y_m/S}) = 0$,
we obtain an exact sequence
$$
\align
&0 \to f_{m*}\Cal O_{Y_m}(\ulcorner kP_m/m \urcorner + K_{Y_m/S} - Y'_{m,0})
\to f_{m*}\Cal O_{Y_m}(\ulcorner kP_m/m \urcorner + K_{Y_m/S}) \\
&\to \bar f_{m,0*}(\Cal O_{Y_m}(\ulcorner kP_m/m \urcorner + K_{Y_m/S})
\otimes \Cal O_{Y'_{m,0}}) \to 0.
\endalign
$$
On the other hand, since $\Cal O_{Y_m}(K_{Y_m/S}) \otimes \Cal O_{Y'_{m,0}}
= \omega_{Y'_{m,0}}(f_m^*E)$, 
there is a natural injective homomorphism
$$
\align
&f_{m,0*}\Cal O_{X'_{m,0}}(\ulcorner kP_m/m \vert_{X'_{m,0}} \urcorner 
+ K_{X'_{m,0}}) \\
&\to \bar f_{m,0*}(\Cal O_{Y_m}(\ulcorner kP_m/m \urcorner + K_{Y_m/S})
\otimes \Cal O_{Y'_{m,0}})(- E).
\endalign
$$
Therefore, there is an injective homomorphism
$$
\Cal I_{kM}((k+1)K_{X'_0}) \to 
\Cal I_{k\bar M}((k+1)K_{Y/S}) \otimes \Cal O_{X'_0}(- E).
$$
Since $R^1\mu_*(\Cal I_{k\bar M}((k+1)K_{Y/S})) = 0$,
the homomorphism
$$
\mu_*(\Cal I_{k\bar M}((k+1)K_{Y/S})) 
\to \mu_*(\Cal I_{k\bar M}((k+1)K_{Y/S}) \otimes \Cal O_{Y_0})
$$
is surjective.
On the other hand, 
since the natural homomorphism $\mu_*\Cal I_{M_{k+1}}((k+1)K_{Y/S}) \to
\mu_*\Cal O_Y((k+1)K_{Y/S})$ is bijective, so is 
$\mu_*\Cal I_{k\bar M}((k+1)K_{Y/S}) \to \mu_*\Cal O_Y((k+1)K_{Y/S})$.
Since we have an injective homomorphism
$\Cal O_{X'_0}(-E) \to \Cal O_{Y_0}$, the lemma is proved.
\hfill $\square$
\enddemo

\proclaim{Lemma 5}
Assume that there exists a divisor $C$ on $X'_0$ such that 
$\Cal I_{kN}(- C) \subset \Cal I_{kM}$ for all 
positive integers $k$. Then the natural homomorphism
$\mu_{0*}(\Cal I_{kM}(kK_{X'_0})) \to \mu_{0*}\Cal O_{X'_0}(kK_{X'_0})$ 
is bijective for any positive integer $k$.
\endproclaim

\demo{Proof}
If the lemma does not hold, then there exists a positive integer $k$ such that
some of the coefficients of 
$$
g_{km,0}^*N_k - kM_{m,0}/m - \Delta_{m,0} + F_{m,0}
$$
at the irreducible components of $F_{m,0}$ are not positive
for any positive integer $m$ with $m \ge k$.
Indeed, if 
$\ulcorner - kM_{m,0}/m \urcorner - \Delta_{m,0} \ge - g_{km,0}^*N_k$,
then the homomorphism
$\mu_{0*}(\Cal I_{kM}(kK_{X'_0})) \to \mu_{0*}\Cal O_{X'_0}(kK_{X'_0})$
is surjective.

Since $\Cal I_{k'N}(- C) \subset \Cal I_{k'M}$ for any $k'$,
there exist positive integers $m'$, depending on $k'$, such that 
$$
g_{k'm',0}^*N_{k'} + f_{m',0}^*C - k'M_{m',0}/m' - \Delta_{m',0} + F_{m',0}
$$
has strictly positive coefficients at the irreducible components of $F_{m',0}$,
because $\Cal O_{X'_{k',0}}(- N_{k'}) \subset \Cal I_{k'N}\Cal O_{X'_{k',0}}$.
Since $g^*_{kk',0}N_k/k \ge N_{k'}/k'$ if $k \vert k'$, the coefficients of
$$
\align
&g_{km',0}^*N_k - kM_{m',0}/m' - \Delta_{m',0} + F_{m',0} \\
&+ \frac{k}{k'}(f_{m',0}^*C + \frac{k'-k}k\Delta_{m',0} - \frac{k'-k}kF_{m',0})
\endalign
$$
are positive.

We may assume that the support of $f_{k,0}^*C$ is contained in $F_{k,0}$.
Then
$$
- f_{k,0}^*C - \frac{k'-k}k\Delta_{k,0} + \frac{k'-k}kF_{k,0} 
$$
is effective for sufficiently large $k'$. 
Since $g_{km',0*}(- \Delta_{k,0} + F_{k,0}) 
\le - \Delta_{m',0} + F_{m',0}$ for $k \le m'$, it follows that
$$
- f_{m',0}^*C - \frac{k'-k}k\Delta_{m',0} + \frac{k'-k}kF_{m',0} 
$$
is also effective for such $k'$ and for any $m' \ge k$.
But this is a contradiction.
\hfill $\square$
\enddemo

\demo{Proof of Main Theorem}
We have to prove that the natural homomorphism
$$
\mu_*\Cal O_Y(mK_{Y/S}) \to \mu_{0*}\Cal O_{X'_0}(mK_{X'_0})
$$
is surjective for any positive integer $m$.
Suppose the contrary. Then there exists a positive integer $k$ 
such that 
$\mu_{0*}(\Cal I_{kM}(kK_{X'_0})) \ne \mu_{0*}\Cal O_{X'_0}(kK_{X'_0})$ 
by Lemma 4.

We consider the conditions
$$
\Cal I_{(k-1)N}(- B_0) \subset \Cal I_{(k+m_0-1)M}
$$
for positive integers $k$.
For $k = 1$, the condition is satisfied, because
$\Cal O_{X'_0}(- B_0) \subset \Cal I_{m_0M}$.
By Lemma 5, there exists a positive integer $k$ such that
$\Cal I_{(k-1)N}(- B_0) \subset \Cal I_{(k+m_0-1)M}$ but
$\Cal I_{kN}(- B_0) \not\subset \Cal I_{(k+m_0)M}$.

We take $A$ such that $A_0 - K_{X'_0}$ for $A_0 = A \vert_{X'_0}$ 
satisfies Theorem 2 for $\mu_0: X'_0 \to X_0$.
Since $A_0 + B_0 \sim m_0K_{X'_0}$, we have
$\Cal I_{kN}(- B_0 + (k+m_0)K_{X'_0})
\simeq \Cal I_{kN}(A_0 - K_{X'_0} + (k+1)K_{X'_0})$.
Thus the homomorphism 
$$
\mu_0^*\mu_{0*}(\Cal I_{kN}(-B_0+(k+m_0)K_{X'_0})) 
\to \Cal I_{kN}(-B_0+(k+m_0)K_{X'_0})
$$
is surjective by Theorem 2.

By Lemma 4, we have
$$
\align
&\mu_{0*}(\Cal I_{kN}(-B_0+(k+m_0)K_{X'_0})) \\
&\subset \mu_{0*}(\Cal I_{(k-1)N}(-B_0+(k+m_0)K_{X'_0})) \\
&\subset \mu_{0*}(\Cal I_{(k+m_0-1)M}((k+m_0)K_{X'_0})) \\
&= \mu_{0*}(\Cal I_{(k+m_0)M}((k+m_0)K_{X'_0})). 
\endalign
$$
But 
$$
\Cal I_{kN}(-B_0+(k+m_0)K_{X'_0})
\not\subset \Cal I_{(k+m_0)M}((k+m_0)K_{X'_0}),
$$
a contradiction.
\hfill $\square$
\enddemo

The following is a slight generalization of a theorem of Siu [Si]:

\proclaim{Theorem 6}
Let $\pi: X \to S$ be a projective flat morphism from a normal variety 
to a germ of a smooth curve $(S, s_0)$.
Assume that the fibers $X_s = \pi^{-1}(s)$
have only canonical singularities and are of general type for all $s \in S$.
Then the plurigenus $P_m(X_s) = \text{dim }H^0(X_s, mK_{X_s})$ is constant
as a function on $s \in S$ for any positive integer $m$.
\endproclaim

\demo{Proof}
The proof is very similar to that of Main Theorem.
We only point out how to modify the proof.
First, the morphisms $f_m: Y_m \to Y$ are constructed only for those $m$ 
such that $P_m(X_{s_0}) \ne 0$.  Moreover, the conditions (1) through (7) are 
assumed only if $P_m(X_s) \ne 0$ for general fibers $X_s$.
In the conditions (2) and (9), $P_m$ and $Q_m$ are 
$\pi_m$-free and free instead of $f'_m$-free and $f'_{m,0}$-free, respectively.
In the conditions (4) and (11), $M_m$ and $N_m$ are the
$\pi_m$-fixed part and the fixed part 
instead of the $f'_m$-fixed part and the $f'_{m,0}$-fixed part, respectively.
The divisor $A$ is very $\pi'$-ample instead of very $\mu$-ample.
The reason for $\text{Supp }B \not\supset X'_0$ is the following: 
since $- \Gamma$ is effective, it follows otherwise that $B \ge Y_0$, 
and we can replace $B$ by a smaller divisor.
The homomorphism in the formula before Definition 1 is replaced 
by the homomorphism $\pi'_*\Cal O_Y(mK_{Y/S}) \to H^0(X'_0, mK_{X'_0})$.

In the definitions of $\Cal I_{kN}$, $\Cal I_{kM}$ and $\Cal I_{k\bar M}$,
the union is taken over those $m$ which we considered.
In Lemmas 4 and 5, their proofs and in the proof of Main Theorem, 
$\mu_*$ and $\mu_{0*}$ are replaced by $\pi'_*$ and $H^0$.
In the proof of Lemma 5, we take $k$ such that $P_k(X_{s_0}) \ne 0$, and
we assume that $k \vert k'$.
In the proof of Main Theorem, the morphism $\mu_0$ is replaced by the 
structure morphism $X'_0 \to \text{Spec }\Bbb C$.
The homomorphism in the formula in the middle of the proof of
Main Theorem is replaced by the statement that the sheaf
$\Cal I_{kN}(- B_0 + (k + m_0)K_{X'_0})$ is generated by global sections.
\hfill $\square$
\enddemo

\enddocument